\documentclass[preprint,12pt,3p]{elsarticle}

\usepackage{amssymb}
\usepackage{amsmath,amsthm}
\usepackage{mathrsfs}
\usepackage{hyperref}
\usepackage{cleveref}
\usepackage[all,cmtip]{xy}
\usepackage{epsf, graphicx}
\usepackage{latexsym,amsfonts,amsbsy,amssymb}
\usepackage{geometry}
\usepackage{titletoc}
\usepackage{rotating}
\usepackage{algorithm}
\usepackage{algorithmic}
\usepackage{amscd}

\makeatletter

\@addtoreset{equation}{section} \makeatother
\newtheorem{theorem}{Theorem}[section]

\newtheorem{conjecture}[theorem]{Conjecture}
\newtheorem{corollary}[theorem]{Corollary}
\newtheorem{remark}[theorem]{Remark}

\newtheorem{example}[theorem]{Example}
\newtheorem{proposition}[theorem]{Proposition}
\newtheorem{question}[theorem]{Question}

\def\Irr{{\rm Irr}}

\def\O{\mathcal{O}}

\def\F{\mathbb{F}}

\def\Z{\mathbb{Z}}

\def\GL{{\rm GL}}

\def\Q{\mathbb{Q}}
\def\G{\mathbf{G}}
\def\L{\mathbf{L}}
\def\P{\mathbf{P}}
\def\V{\mathbf{V}}
\def\Y{\mathbf{Y}}

\def\GG{{\rm G}\Gamma}
\def\RG{R\Gamma}

\makeatletter
\def\ps@pprintTitle{%
\let\@oddhead\@empty
\let\@evenhead\@empty
\def\@oddfoot{\reset@font\hfil\thepage\hfil}
\let\@evenfoot\@oddfoot
}
\makeatother

\begin{document}

\begin{frontmatter}

\title{Descent of splendid Rickard equivalences in ${\rm GL}_n(q)$}

\author[label1,label2]{Xin Huang}

\address[label1]{School of Mathematics and Statistics, Central China Normal University, Wuhan 430079, China}
\address[label2]{SICM, Southern University of Science and Technology, Shenzhen 518055, China}
\address[label3]{Yau Mathematical Sciences Center, Tsinghua University, Beijing 100084, China}
\address[label4]{School of Mathematical Sciences, Peking University, Beijing 100871, China}

\ead{xinhuang@mails.ccnu.edu.cn}

\author[label3]{Pengcheng Li}
\ead{pcli@tsinghua.edu.cn}

\author[label2,label4]{Jiping Zhang}
\ead{jzhang@pku.edu.cn}


\begin{abstract}

Let $n$ be a positive integer and $q$ a prime power. We prove that a refined version of Brou\'{e}'s abelian defect group conjecture holds for unipotent $\ell$-blocks of ${\rm GL}_n(q)$, where $\ell\nmid q$. We also give a sufficient condition on general $\ell$-blocks of ${\rm GL}_n(q)$ to satisfy the refined abelian defect group conjecture. We explain by an example that
this sufficient condition does not hold in general.

\end{abstract}


\end{frontmatter}


\section{Introduction}\label{s1}

In \cite{Kessar_Linckelmann}, Kessar and Linckelmann proposed a refined version of Brou\'{e}'s abelian defect group conjecture.

\begin{conjecture}[The refined Brou\'{e} conjecture]\label{KLconj}
For an arbitrary complete discrete valuation ring $\O$ and a block $b$ of a finite group $G$ over $\O$ with an abelian defect group, there is a splendid Rickard equivalence between $\O Gb$ and its Brauer correspondent.
\end{conjecture}

Conjecture \ref{KLconj} extends Brou\'{e}'s conjecture because the original conjecture is with the assumption that the complete discrete valuation rings have splitting residue fields.  For blocks with abelian defect groups, Kessar and Linckelmann (\cite[Corollary 1.9]{Kessar_Linckelmann}) showed that Conjecture \ref{KLconj} implies Navarro's refinement of the Alperin--McKay conjecture (\cite[Conjecture B]{Navarro}). This implication has been generalised by Boltje (see \cite[Theorem 1.4]{Boltje}), who proved that Conjecture \ref{KLconj} implies Turull's refinement of the Alperin--McKay conjecture (\cite[Conjecture]{T13}). Note that Turull's refinement of the Alperin--McKay conjecture contains the refined versions of the Alperin--McKay conjecture proposed
by Isaacs--Navarro (\cite[Conjecture B]{IN}) and Navarro (\cite[Conjecture B]{Navarro}). Moreover, by \cite[Theorem 1]{H24}, for blocks with abelian defect groups, Conjecture \ref{KLconj} implies Navarro's refinement of Alperin's weight conjecture (see \cite[Conjecture 2]{H24}). The following alternative version of Conjecture \ref{KLconj} looks slightly weaker, but in fact it is equivalent to Conjecture \ref{KLconj}.

\noindent\textbf{Conjecture 1.1$'$.} {\it For an arbitrary complete discrete valuation ring $\O$ of characteristic $0$ with residue field of characteristic $\ell$ and a block $b$ of a finite group $G$ over $\O$ with an abelian defect group, there is a splendid Rickard equivalence between $\O Gb$ and its Brauer correspondent.}

\medskip Clearly Conjecture \ref{KLconj} implies Conjecture 1.1$'$. Assume now that Conjecture 1.1$'$ holds. Let $\O$ be a complete discrete valuation ring of characteristic $\ell$ with residue field $k$. Let $G$ be a finite group, $b$ a block of $\O G$ with an abelian defect group $P$, and $c$ the block of $\O N_G(P)$ corresponding to $b$ via the Brauer correspondence. We need to show that there is a splendid Rickard equivalence between $\O Gb$ and $\O N_G(P)c$.  The following argument is inspired by \cite[6.11]{WZZ}. Without loss of generality we may assume that $k$ is perfect. By \cite[Chapter II, Theorem 3]{Serre}, there is a complete discrete valuation ring $\tilde{\O}$ with characteristic $0$ with residue field $k$.
Denote by $\bar{b}$ the image of $b$ in $\O G$ and $\bar{c}$ the image of $c$ in $kN_G(P)$. The blocks $\bar{b}$ and $\bar{c}$ can be lifted to blocks $\tilde{b}$ and $\tilde{c}$ of $\tilde{\O}G$ and $\tilde{\O}N_G(P)$ respectively. By our assumption, there is a splendid Rickard equivalence between $\tilde{\O} G\tilde{b}$ and $\tilde{\O}N_G(P)\tilde{c}$. Hence there is a splendid Rickard equivalence between $kG\bar{b}$ and $\bar{\O}N_G(P)\bar{c}$. Since both $\O$ and $k$ have the same characteristic $\ell$, by \cite[Chapter II, Proposition 8]{Serre}, $k$ can be identified with a subring of $\O$. Hence $\O Gb\cong \O\otimes_k kG\bar{b}$ and $\O N_G(P)c\cong \O\otimes_k kN_G(P)\bar{c}$. So there is a splendid Rickard equivalence between $\O Gb$ and $\O N_G(P)c$ (see e.g. \cite[Proposition 4.5]{Kessar_Linckelmann}).

In this paper we investigate Conjecture \ref{KLconj} for finite general linear groups. We choose to study this family of groups because splendid Rickard equivalences for unipotent blocks of $\GL_n(q)$ constructed by Chuang and Rouquier (\cite[Theorem 7.18]{CR}) are obtained by using the $\mathfrak{sl}_2$-categorification. Splendid Rickard equivalences for blocks of symmetric groups are also obtained by using the $\mathfrak{sl}_2$-categorification, but in that case, any field is a splitting field for symmetric groups, hence Conjecture \ref{KLconj} is already proved by Chuang and Rouquier (\cite[Theorem 7.6]{CR}). 

Let us fix some notation. Throughout the rest of this paper $\ell$ is a prime, $n$ is a positive integer, $q$ is a prime power, $k\subseteq k'$ are fields of characteristic $\ell$, and $\O\subseteq\O'$ are complete discrete valuation rings of characteristic $0$ with $J(\O)\subseteq J(\O')$ and with residue fields $k$ and $k'$, respectively. Denote by $K'$ the quotient field of $\O'$. Assume that $K'$ contains a primitive $|G|$-th root of unity for every finite group $G$ considered below. So both $K'$ and $k'$ are splitting fields for all finite groups considered below. Denote by $\Q$ and $\F_\ell$ the prime fields of $K'$ and $k'$, respectively. Let $\Z_\ell$ be the ring of $\ell$-adic integers. By \cite[Chapter 2, Theorem 3, 4 and Proposition 1]{Serre}, $\Z_\ell$ can be identified with the unique complete discrete valuation ring $R$ contained in $\O$ such that $J(R)=\ell R$ and the image of $R$ under the canonical surjection $\O\to k$ is $\F_\ell$. We always make this identification.

Let $G$ be a finite group. By a {\it block} of the a finite group algebra $\Lambda G$, where $\Lambda\in \{\O,k\}$, we usually mean a primitive idempotent $b$ of the center of $\Lambda G$, and $\Lambda Gb$ is called a {\it block algebra}. Sometimes the term ``block" will also be used to mean the correspondent set of irreducible characters. For a subgroup $H$ of $G$, let $(\Lambda Gb)^H$ denote the set of $H$-fixed elements of the block algebra $\Lambda Gb$ under the conjugation action. If $H$ is a $p$-subgroup, the {\it Brauer map} is the $\Lambda$-algebra homomorphism
${\rm Br}_H: (\Lambda Gb)^H\to kC_G(H)$, $\sum_{g\in G}\alpha_gg\mapsto \sum_{g\in C_G(H)}\bar{\alpha}_gg,$
where $\bar{\alpha}_g$ denotes the image of $\alpha_g$ in $k$.
For a block $b$ of $\Lambda G$, a {\it defect group} of $b$ is a maximal $p$-subgroup $P$ of $G$ such that ${\rm Br}_P(b)\neq 0$. By Brauer's first main theorem, there is a unique block $c$ of $\Lambda N_G(P)$ with defect group $P$ such that ${\rm Br}_P(b) = {\rm Br}_P(c)$ and the map $b\mapsto c$ is a bijection between the set of blocks of $\Lambda G$ with defect group $P$ and the set of blocks of $\Lambda N_G(P)$ with defect group $P$. This bijection is known as the {\it Brauer correspondence}.

Let $b$ be a block $\O \GL_n(q)$ with a defect group $P$. Let $b'$ be a block of $\O'\GL_n(q)$ with $bb'=b'=b'b$, then it is easy to show that $P$ is also a defect group of $b'$. Assume that $P$ is abelian. If $\ell\mid q$, then by the last paragraph of \cite[Section 4]{Dagger}, $P$ is either the trivial subgroup or a Sylow $\ell$-subgroup of $\GL_n(q)$. So if $P\neq 1$, then $n=1$ or $n=2$ (otherwise, the Sylow $\ell$-subgroups are not abelian). If $n=1$ and $\ell\mid q$, Conjecture \ref{KLconj} is trivially true for the group $\GL_n(q)$. If $n=2$ and $\ell\mid q$, Conjecture \ref{KLconj} is also true for $\GL_n(q)$ according to \cite[Theorem 1.1]{HLZ}.
So for the group $\GL_n(q)$ we only need to consider Conjecture \ref{KLconj} in non-defining characteristic. From now on, we assume that $\ell\nmid q$.

For the large enough field $k'$, Chuang and Rouquier (\cite[Theorem 7.20]{CR}) proved that there is a splendid Rickard equivalence between every block algebra of $k'\GL_n(q)$ with an abelian defect group and its Brauer correspondent algebra. The main result of this paper is the following.

\begin{theorem}\label{main}
Let $G$ be the group $\GL_n(q)$, $b$ a unipotent block of $\O'G$ with an abelian defect group $P$, and $c$ the block of $\O'N_G(P)$ corresponding to $b$ via the Brauer correspondence. Then $b\in \Z_\ell G$, $c\in \Z_\ell N_G(P)$, and the block algebras $\Z_\ell Gb$ and $\Z_\ell N_G(P)c$ are splendidly Rickard equivalent. More precisely, there is a splendid Rickard complex $X$ of $(\Z_\ell Gb,\Z_\ell N_G(P)c)$-bimodules such that $\O'\otimes_{\Z_\ell}X$ is isomorphic to Chuang and Rouquier's complex $X'$.
\end{theorem}

See \S\ref{splendid} below for the definition of a splendid Rickard equivalence. By \cite[Theorem 1.4]{Boltje} and \cite[Theorem 1]{H24}, Theorem \ref{main} has the following corollary.

\begin{corollary}
\begin{enumerate}[\rm (i)]
	\item Turull's refinement of the Alperin--McKay conjecture (\cite[Conjecture]{T13}) holds for unipotent $\ell$-blocks of $\GL_n(q)$ with an abelian defect group.
	\item The block version of Navarro's refinement of Alperin's weight conjecture (\cite[Conjecture 2]{H24}) holds for unipotent $\ell$-blocks of $\GL_n(q)$ with an abelian defect group.
\end{enumerate}

\end{corollary}

By lifting theorem of splendid Rickard equivalences (reviewed in Theorem \ref{lifting} below), to prove Theorem \ref{main}, we may replace $\O'$ by $k'$ and replace $\Z_\ell$ by $\F_\ell$. The main steps of proving Theorem \ref{main} are summarized in the following remarks.

\begin{remark}\label{1.4}
{\rm Let $G$ be $\GL_n(q)$ and $b$ a unipotent block of $k'G$. Let $e$ be the multiplicative order of $q$ in $k^\times$. Recall that there is associated a non-negative integer $w$, called the ($e$-)weight of $b$. Let $P$ be a defect group of $b$, and let $c$ be the block of $k'N_G(P)$ corresponding to $b$ via the Brauer correspondence. The defect group $P$ of $b$ is abelian if and only if $w<\ell$. The proof of \cite[Theorem 7.20]{CR} consists of three steps:

\noindent (i). $k'N_G(P)c$ is splendidly Rickard equivalent to the principal block algebra of $k'(\GL_e(q)\wr S_w)$ (by \cite[Theorem 10.1]{Rou} and \cite[Theorem 4.3 (b)]{Marcusonequivalences} ).

\noindent (ii). There exists an integer $n'\geq1$ and a unipotent block $b'$ of $\GL_{n'}(q)$, with weight $w$, which is splendidly Morita equivalent to the principal block algebra of $k'(\GL_e(q)\wr S_w)$ (see \cite[Theorem 5.0.7]{Mi01} or \cite[Theorem 1]{Turner}).

\noindent (iii). The block algebras $k'\GL_{n'}(q)b'$ and $k'Gb$ are splendidly Rickard equivalent (see \cite[Theorem 7.18]{CR}). }
\end{remark}

\begin{remark}\label{1.5}
{\rm Assume that $b$ is a unipotent block of $k'G$. By Corollary \ref{unipotent} below, we have $b\in \F_\ell G$. By Remark \ref{1.4}, there is a splendid Rickard complex $X'_1$ for $k'N_G(P)b$ and $k'\GL_{n'}(q)b'$, and a splendid Rickard complex $X'_2$ for $k'\GL_{n'}(q)b'$ and $k'Gb$. We will show (in Section \ref{s3} and \ref{s4}) that there is a complex $X_1$ of $(\F_\ell N_G(P)b,\F_\ell\GL_{n'}(q)b')$-bimodules, and a complex $X_2$ of $(\F_\ell\GL_{n'}(q)b',\F_\ell Gb)$-bimodules, such that $X'_1\cong k'\otimes_{\F_\ell}X_1$ and $X'_2\cong k'\otimes_{\F_\ell}X_2$.
Then by \cite[Proposition 4.5 (a)]{Kessar_Linckelmann}, $X_1$ induces a splendid Rickard equivalence between $\F_\ell N_G(P)b$ and $\F_\ell \GL_{n'}(q)b'$, and $X_2$ induces a splendid Rickard equivalence between $\F_\ell\GL_{n'}(q)b'$ and $\F_\ell Gb$. Hence at that time, the proof of Theorem \ref{main} is complete.}
\end{remark}

In \cite{DVVcl} and \cite{DVVuni}, Brou\'{e}'s abelian defect group conjecture is proved for unipotent $\ell$-blocks of the groups ${\rm GU}_n(q)$, ${\rm Sp}_{2n}(q)$ and ${\rm SO}_{2n+1}(q)$ at linear primes $\ell$. Since the constructions of the derived equivalences for unipotent blocks of these groups have many common properties with the constructions of the equivalences for unipotent blocks of $\GL_n(q)$, the methods in the proof of Theorem \ref{main} may be used to prove a similar proposition for these groups instead of $\GL_n(q)$.

One will ask whether the refined abelian defect group conjecture holds for general blocks of $\O\GL_{n}(q)$, not only unipotent blocks.  In Section \ref{ap}, we investigate this question and give a sufficient condition under which there exists a splendid Rickard equivalence between a block of $\O\GL_n(q)$ and its Brauer correspondent. We
provide an example to show that this sufficient condition does not hold in general; see Example \ref{example}.

\section{Preliminaries}

\subsection{Notation}

For a finite group $G$, we denote by $G^{\rm op}$ the opposite group and by $\Delta G$ the subgroup $\{(g,g^{-1})~|~g\in G\}$ of $G\times G^{\rm op}$. For an additive category $\mathcal{C}$, we denote by ${\rm Comp}^b(\mathcal{C})$ the category of bounded complexes of objects of $\mathcal{C}$ and by ${\rm Ho}^b(\mathcal{C})$ its homotopy category. For an algebra $A$, we denote by $A^{\rm op}$ the opposite algebra of $A$. Unless specified otherwise, modules in the paper are left modules. We denote by $A$-${\rm mod}$ the category of finitely generated $A$-modules, and by $G_0(A)$ the Grothendieck group of $A$-${\rm mod}$. Let ${\rm Comp}^b(A):={\rm Comp}^b(A\mbox{-mod})$ and ${\rm Ho}^b(A):={\rm Ho}^b(A\mbox{-mod})$.
Let $C\in {\rm Comp}^b(A)$. There is a unique (up to a non-unique isomorphism) complex $C^{\rm red}$ that is isomorphic to $C$ in the homotopy category ${\rm Ho}^b(A)$ and that has no nonzero direct summand that is homotopy equivalent to $0$. So $C\cong C^{\rm red}\oplus C_0$ for some $C_0\in {\rm Comp}^b(A)$ homotopy equivalent to zero (see e.g. \cite[Corollary 1.18.19]{Linckelmann}). 

\subsection{Splendid Rickard equivalences}\label{splendid}

Let $A$ and $B$ be symmetric $\Lambda$-algebras, where $\Lambda\in \{\O,k\}$. Let $X$ be a bounded complex of finitely generated $(A,B)$-bimodules which are projective as left $A$-modules and as right $B$-modules, and let $X^*:={\rm Hom}_\Lambda(X,\Lambda)$ be the dual complex. It is said that $X$ induces a {\it Rickard equivalence} and that $X$
is a {\it Rickard complex} if there exists a contractible complex of $(A,A)$-bimodules $Y$
and a contractible complex of $(B,B)$-bimodules $Z$ such that $X\otimes_B X^*\cong A\oplus Y$ as complexes $(A,A)$-bimodules and $X^*\otimes_A X\cong B\oplus Z$ as complexes of $(B,B)$-bimodules.
Let $G$ and $H$ be finite groups. Let $b$ (resp. $c$) be a block of $\Lambda G$ (resp. $\Lambda H$). Let $X:=(X_n)_{n\in \mathbb{Z}}$ be a Rickard complex of $(\Lambda Gb,\Lambda Hc)$-bimodules. If each $X_n$ is a direct summand of permutation $\Lambda(G\times H^{\rm op})$-module (i.e., an $\ell$-permutation $\Lambda(G\times H^{\rm op})$-module), then $X$ is said to be {\it splendid}; $\Lambda Gb$ and $\Lambda Hc$ are said to be {\it splendidly Rickard equivalent}. Chuang and Rouquier remarked in the first paragraph of \cite[\S7.1.2]{CR} that one usually puts some condition on the vertex of $X_n$, but this is actually automatic. The following theorem on lifting splendid Rickard equivalences is due to Rickard.

\begin{theorem}[{\cite[Theorem 5.2]{Rickardsplendid}}]\label{lifting}
Let $G$ and $H$ be finite groups. Let $b$ (resp. $c$) be an idempotent in the center of $\O G$ (resp. $\O H$). Denote by $\bar{b}$ (resp. $\bar{c}$) the image of $b$ (resp. $c$) in $kG$ (resp. $kH$). Assume that there is a complex $\bar{X}$ of $(kG\bar{b}, kH\bar{c})$-bimodules inducing a splendid Rickard equivalence. Then there is a complex $X$ of $(\O Gb, \O Hc)$-bimodules inducing a splendid Rickard equivalence and satisfying $k\otimes_\O X\cong \bar{X}$.
\end{theorem}

Note that although the statement in \cite[Theorem 5.2]{Rickardsplendid} is for principal blocks, but
the proof carries over nearly verbatim to arbitrary blocks. We also note that the blanket
assumption in [12] that the coefficient rings are big enough is not used in the proof of \cite[Theorem 5.2]{Rickardsplendid}.

\subsection{Actions of Galois automorphisms on modules}
Let $R\subseteq R'$ be two commutative domains. Let $A$ be an $R$-algebra and let $A':=R'\otimes_R A$. Let $\Gamma$ be the group of automorphisms of $R'$ which restricts to the identity map on $R$.  For an $A'$-module
$U'$ and a $\sigma\in \Gamma$, denote by ${}^\sigma U'$ the $A'$-module which is equal to $U'$ as a module over the subalgebra $1\otimes A$ of $A'$, such that $x\otimes a$ acts on $U'$ as $\sigma^{-1}(x) \otimes a$ for all $a \in A$ and $x\in R'$.
 The $A'$-module $U'$ is {\it $\Gamma$-stable} if ${}^\sigma U'\cong U'$ for all $\sigma\in \Gamma$. $U'$ is said to be {\it defined over $R$}, if there is an $A$-module $U$ such that $U'\cong R'\otimes_R U$.
In this special case, $U'$ is $\Gamma$-stable, because for any $\sigma\in\Gamma$, the map sending $x\otimes u$ to $\sigma^{-1}(x)\otimes u$ is an isomorphism $R'\otimes_R U\cong {}^\sigma(R'\otimes_R U)$, where $u\in U$ and $x\in R'$. 

\subsection{Block idempotents and coefficient rings}

The following result can be deduced by \cite[Lemma 2.2 (b)]{Farrell}. For the convenience of the reader we include a proof.

\begin{proposition}\label{2.2}
Let $G$ be a finite group, $b$ a block of $\O'G$, and let $\chi:G\to K'$ be the character of a simple $K'Gb$-module. Let $\bar{b}$ be the image of $b$ in $k'G$. If the values of $\chi$ are contained in $\Q$ (hence in $\mathbb{Z}$), then we have $b\in\Z_\ell G$ and $\bar{b}\in \F_\ell G$.
\end{proposition}

\noindent{\it Proof.} Let $V$ be an $\O'G$-module such that the $K'G$-module $K'\otimes_{\O'}V$ affords the character $\chi$ (see e.g. \cite[Theorem 4.16.5]{Linckelmann} for the existence of $V$). Let $\varphi:G\to k'$ be the character afforded by the $k'G\bar{b}$-module $k'\otimes_{\O'}V$. The values of $\varphi$ are the images of values of $\chi$ under the canonical surjection $\O\to k$, hence contained in $\F_\ell$. For any $\sigma\in {\rm Gal}(k'/\F_\ell)$, $\sigma$ induces a ring automorphism of $k'G$ in an obvious way. Hence $\sigma(\bar{b})$ is also a block of $k'G$. Since $\varphi$ is invariant under the action of $\sigma$, we see that $\sigma(\bar{b})=\bar{b}$. Since every finite group has a finite splitting field, we may assume that $k'$ is finite. Then we deduce that $\bar{b}\in \F_\ell G$ by the Galois theory. By idempotent lifting arguments, we have $b\in \Z_\ell G$. $\hfill\square$

\begin{corollary}\label{unipotent}
Let $b$ be a unipotent block of $\O'\GL_n(q)$, and let $\bar{b}$ be the image of $b$ in $k'\GL_n(q)$, then $b\in \Z_\ell \GL_n(q)$ and $\bar{b}\in \F_\ell \GL_n(q)$.
\end{corollary}

\noindent{\it Proof.} Let $\chi:\GL_n(q)\to K'$ be a unipotent character of $\GL_n(q)$. By \cite[Example 1.1]{Geck}, the values of $\chi$ are contained in $\Q$. The statement follows by Proposition \ref{2.2}. $\hfill\square$

\medskip It is easy to see that Proposition \ref{2.2} also has the following corollary.

\begin{corollary}\label{principal}
Let $G$ be a finite group and let $b$ be the principal block of $\O'G$, then $b\in\Z_\ell G$ and $\bar{b}\in \F_\ell G$.
\end{corollary}

\section{On unipotent blocks of general linear groups with same weights}\label{s3}

Keep the notation of Remark \ref{1.4} and \ref{1.5}. In this section, we show that the splendid Rickard complex $X'_2$ is defined over $\F_\ell$. The construction of this complex is played back to \cite[Theorem 7.18]{CR}.

We start with the case $\ell|(q-1)$. By \cite[Remark 7.19]{CR}, $k'G$ ($=k'\GL_n(q)$) has a unique unipotent block $b$, the principal block. The number of isomorphism classes of simple $k'Gb$-modules is the number of partitions of $n$. So if $n\neq n'$, $k'Gb$ can not be Rickard equivalent to a unipotent block of $k'\GL_{n'}(q)$.
Hence $n=n'$. In this case, the complex $X'_2$ can be taken to be the $(k'Gb,k'Gb)$-bimodule $k'Gb$, and it is obviously defined over $\F_\ell$.

It remains to consider the case where $\ell\nmid q(q-1)$. The construction of the complex $X_2'$ uses the $\mathfrak{sl}_2$-categorification. Let $G_n:=\GL_n(q)$ and let $A'_n=k'G_nb_n$ be the sum of the unipotent block algebras of $k'G_n$. Given a finite group $H$ with $\ell\nmid|H|$, put $e_H:=\frac{1}{|H|}\sum_{h\in H}h$, an idempotent in $\F_\ell H\subseteq k'H$. For a matrix $g\in G_n$, denote by ${}^tg$ the transpose of $g$. Denote by $V_n$ the subgroup of upper triangular matrices of $G_n$ with diagonal coefficients $1$ whose off-diagonal coefficients vanish outside the $n$-th column. Denote by $D_n$ the subgroup of $G_n$ of diagonal matrices with diagonal entries $1$ except the $(n,n)$-th one.
\[V_n=\left( {\begin{array}{*{20}{c}}
  1& & &* \\
   & \ddots & & \vdots  \\
   & &1&* \\
   & & & 1
\end{array}} \right)_{i\times i},~~~~~~ D_n=\left( {\begin{array}{*{20}{c}}
  1& & & \\
   & \ddots & &   \\
   & &1& \\
   & & & *
\end{array}} \right)_{i\times i}.\]

For $i\in\{0,1,\cdots, n-1\}$, we view $G_i$ as a subgroup of $G_n$ via the first $i$ coordinates. Following Chuang and Rouquier, we put
$$E'_{i,n}:=k'G_ne_{(V_n\rtimes\cdots\rtimes V_{i+1})\rtimes(D_{i+1}\times\cdots\times D_n)}\otimes_{k'G_i}-:A'_i\mbox{-}{\rm mod}\to A'_n\mbox{-}{\rm mod}$$
and put
$$F'_{i,n}:=e_{(V_n\rtimes\cdots\rtimes V_{i+1})\rtimes(D_{i+1}\times\cdots\times D_n)}k'G_n\otimes_{k'G_n}-:A'_n\mbox{-}{\rm mod}\to A'_i\mbox{-}{\rm mod}.$$
The functors $E'_{i,n}$ and $F'_{i,n}$ are canonically left and right adjoint.
Let $\mathcal{A}':=\bigoplus_{n\geq0}A'_n$-mod, $E':=\bigoplus_{n\geq0}E'_{n,n+1}$ and $F':=\bigoplus_{n\geq0}F'_{n,n+1}$. Denote by $X$ the endomorphism of $E'$ given on $E'_{n-1,n}$ by right multiplication by
$$\hat{X}_n:=q^{n-1}e_{V_nD_n}e_{{}^tV_n}e_{V_nD_n}.$$

Given $a\in k'^{\times}$, let $E'_a$ be the generalised $a$-eigenspace of $X$ on $E'$. We have a decomposition $E'=\bigoplus_{a\in k'^\times}E'_a$. There is a corresponding decomposition $F'=\bigoplus_{a\in k'^\times}F'_a$, such that $F'_a$ is left and right adjoint to $E'_a$. Note that $E'_a$ and $F'_a$ are functors from $\mathcal{A}'$ to $\mathcal{A}'$, so they induce actions $[E'_a]$ and $[F'_a]$ on $\bigoplus_{n\geq 0}G_0(A'_n\mbox{-}{\rm mod})$, respectively.  By \cite[Lemma 7.16]{CR}, the action of $[E'_a]$ and $[F'_a]$ on $\bigoplus_{n\geq 0}G_0(A'_n\mbox{-}{\rm mod})$ gives a representation of $\mathfrak{sl}_2$, and the classes of simple objects are weight vectors. Moreover, these actions give rise to an $\mathfrak{sl}_2$-categorification on $\mathcal{A}'$ (see \cite[\S7.3.1]{CR}).

The functors $E'$ and $F'$ are defined by tensoring with bimodules. Since bimodules are more convenient to handle than functors, let us add some notation. Let $\mathscr{A}':=\bigoplus_{n\geq0}A'_n$, a $k'$-algebra. The functor $E'_{n-1,n}$ is defined by the $(k'G_n,k'G_{n-1})$-bimodule $k'G_ne_{V_nD_n}$, and we denote this bimodule by $\mathscr{E}'_{n-1,n}$. We can view $\mathscr{E}'_{n-1,n}$ as an $(\mathscr{A}',\mathscr{A}')$-bimodule, so we have an $(\mathscr{A}',\mathscr{A}')$-bimodule
$$\mathscr{E}':=\bigoplus_{n\geq0}k'G_ne_{V_nD_n}.$$
Clearly the bimodule $\mathscr{E}'$ corresponds to the functor $E'$. Similarly, we have an $(\mathscr{A}',\mathscr{A}')$-bimodule
$$\mathscr{F'}:=\bigoplus_{n\geq0}e_{V_nD_n}k'G_n,$$
which corresponds to the functor $F'$.
The endomorphism $X$ of $E'$ has similar properties with the endomorphism $\mathscr{X}$ of $\mathscr{E}'$, given on $\mathscr{E}'_{n-1,n}$ by right multiplication by
$$\hat{X}_n=q^{n-1}e_{V_nD_n}e_{{}^tV_n}e_{V_nD_n}.$$

Given $a\in k'^{\times}$, let $\mathscr{E}'_a$ be the generalised $a$-eigenspace of $\mathscr{X}$ on $\mathscr{E}'$. We have a decomposition $\mathscr{E}'=\bigoplus_{a\in k'^\times}\mathscr{E}'_a$. Then the $(\mathscr{A}',\mathscr{A}')$-bimodule $\mathscr{E}'_a$ corresponds to the functor $E'_a$. Since the bimodule $\mathscr{F}'$ corresponds to the functor $F'$, there is a decomposition $\mathscr{F}'=\bigoplus_{a\in k'^\times}\mathscr{F}'_a$, such that the $(\mathscr{A}',\mathscr{A}')$-bimodule $\mathscr{F}'_a$ corresponds to the functor $F'_a$.

\begin{proposition}\label{eigenvalue}
The eigenvalues of $X$ as an endomorphism of $E'$ are contained in $\F_\ell$. Hence the eigenvalues of $\mathscr{X}$ as an endomorphism of $\mathscr{E}'$ are contained in $\F_\ell$.
\end{proposition}

\noindent{\it Proof.} The set of eigenvalues of $X$ is the union of the sets of eigenvalues of $\hat{X}_n$'s. Note that $\hat{X}_1$ induces the identical map on $E'_{1,2}$, so each eigenvalue of $\hat{X}_1$ is $1$. By \cite[Lemma 4.7]{Gr}, the eigenvalues of $\hat{X}_n$ on $E'_{n,n-1}$ are powers of $q$ (considered as elements in $\F_\ell$), whence the statement. $\hfill\square$

\medskip Since the results on the local block theory of symmetric groups generalise to unipotent blocks of general linear groups  \cite[\S3]{Broue86}, we have an analog of Theorem 7.1 in \cite{CR}, which is omitted in \cite{CR}.

\begin{theorem}[{an analog of \cite[Theorem 7.1]{CR}}]\label{LLT}
Let $e$ be the multiplicative order of $q$ in $\F_\ell$. The functors $[E'_a]$ and $[F'_a]$ for $a\in \F_\ell$ give rise to an action of the affine Lie algebra $\hat{\mathfrak{sl}}_e$ on $\bigoplus_{n\geq 0}G_0(A'_n\mbox{-}{\rm mod})$. The decomposition of $\bigoplus_{n\geq 0}G_0(A'_n\mbox{-}{\rm mod})$ in blocks coincides with its decomposition in weight spaces. Two unipotent blocks of general linear groups have the same weight if and only if they are in the same orbit under the adjoint action of the affine Weyl group.
\end{theorem}

In Remark \ref{1.4} and \ref{1.5}, the unipotent block algebras $k'Gb$ and $k'\GL_{n'}(q)b'$ have the same weight. So by Theorem \ref{LLT}, there is a sequence of unipotent block algebras $B'_0=k'Gb, B'_1,\cdots, B'_s=k'\GL_{n'}(q)b'$ such that $B'_j$ is the image of $B'_{j-1}$ by some simple reflection $\sigma_{a_j}$ of the affine Weyl group. By Proposition \ref{eigenvalue}, these eigenvalues $a_1,\cdots,a_s$ are contained in $\F_\ell$.

By \cite[Theorem 6.4]{CR}, the complex of functors $\Theta'$ there associated with $a=a_j$ induces a self-equivalence of ${\rm Ho}^b(\mathcal{A}')$. It restricts to a splendid Rickard equivalence between $B'_j$ and $B'_{j-1}$. The Rickard equivalence between $k'Gb$ and $k'\GL_{n'}(q)b'$ is the composition of these equivalences.
In other words, if we denote by $C'_j$ the complex of $(B'_{j}, B'_{j-1})$-bimodules which induces the splendid Rickard equivalence between $B'_{j}$ and $B'_{j-1}$, then
$$X'_2= C'_s\otimes_{B'_{s-1}}\cdots\otimes_{B'_1}C'_1.$$
By Corollary \ref{unipotent}, there are block algebras $B_0=\F_\ell Gb, B_1,\cdots, B_s=\F_\ell\GL_{n'}(q)b'$ such that $B'_j\cong k'\otimes_{\F_\ell}B_j$.

\begin{proposition}\label{3.3}
There is a complex $X_2$ of $(\F_\ell\GL_{n'}(q)b',\F_\ell Gb)$-bimodules, such that $X'_2\cong k'\otimes_{\F_\ell}X_2$.
\end{proposition}

To prove Proposition \ref{3.3}, it suffices to prove the following statement.

\begin{proposition}
For each $j\in \{1,\cdots, s\}$, there exists a complex $C_j$ of $(B_j,B_{j-1})$-bimodules such that $C'_j\cong k'\otimes_{\F_\ell}C_j$.
\end{proposition}

\noindent{\it Proof.} The categories $B'_{j-1}$-mod and $B'_j$-mod are exactly the categories $\mathcal{A}'_{-\lambda}$ and $\mathcal{A}'_\lambda$ in \cite[Theorem 6.4]{CR} for some $\lambda$. The complex $C'_j$ of $(B'_{j},B'_{j-1})$-bimodules is defined by the complex of functors $\Theta'_\lambda$ described in \cite[\S6.1]{CR}. We can express the complex $C'_j$ in terms of $\Theta'_\lambda$. Since the functor
$$\Theta'_\lambda:{\rm Comp}^b(\mathcal{A}'_{-\lambda})\to {\rm Comp}^b(\mathcal{A}'_{\lambda})$$
coincides with the functor
$$C'_j\otimes_{B'_{j-1}}-:{\rm Comp}^b(B'_{j-1})\to {\rm Comp}^b(B'_j),$$
we have
$$C'_j\cong C'_j\otimes_{B'_{j-1}}B'_{j-1}\cong \Theta'_\lambda(B'_{j-1}).$$
More explicitly, if we write
$$\Theta'_\lambda=\cdots\to(\Theta'_\lambda)^i\xrightarrow{d'^i}(\Theta'_\lambda)^{i+1}\to\cdots,$$
then $C'_j$ is isomorphic to the complex
$$\cdots\to(\Theta'_\lambda)^i(B'_{j-1}) \xrightarrow{d'^i(B'_{j-1})}(\Theta'_\lambda)^{i+1}(B'_{j-1})\to\cdots$$
(see the first paragraph of \cite[\S4.1.4]{CR}). Note that since the functor $(\Theta'_\lambda)^i$ is defined by tensoring with bimodules and since $B'_{j-1}$ is a $(B'_{j-1},B'_{j-1})$-bimodule, $(\Theta'_\lambda)^i(B'_{j-1})$ is not only a left $B'_j$-module but also a right $B'_{j-1}$-module, hence a $(B'_j,B'_{j-1})$-bimodule.

To show that the complex $C'_j$ is defined over $k$, let us first add some notation. By Corollary \ref{unipotent}, we have $b_n\in \F_\ell G_n$ (recall that $b_n$ denotes the sum of the unipotent block of $k'G_n$). Put $A_n:=\F_\ell G_nb_n$, then we have $A'_n\cong k'\otimes_{\F_\ell}A_n$. We put
$$E_{i,n}:=\F_\ell G_ne_{(V_n\rtimes\cdots\rtimes V_{i+1})\rtimes(D_{i+1}\times\cdots\times D_n)}\otimes_{\F_\ell G_i}-:A_i\mbox{-}{\rm mod}\to A_n\mbox{-}{\rm mod}$$
and put
$$F_{i,n}:=e_{(V_n\rtimes\cdots\rtimes V_{i+1})\rtimes(D_{i+1}\times\cdots\times D_n)}\F_\ell G_n\otimes_{\F_\ell G_n}-:A_n\mbox{-}{\rm mod}\to A_i\mbox{-}{\rm mod}.$$
The functors $E_{i,n}$ and $F_{i,n}$ are canonically left and right adjoint. 

Let $\mathcal{A}:=\bigoplus_{n\geq0}A_n$-mod, $E:=\bigoplus_{n\geq0}E_{n,n+1}$ and $F:=\bigoplus_{n\geq0}F_{n,n+1}$. Let $\mathscr{A}:=\bigoplus_{n\geq0}A_n$, a $k$-algebra. The functor $E_{n-1,n}$ is defined by the $(\F_\ell G_n,\F_\ell G_{n-1})$-bimodule $\F_\ell G_ne_{V_nD_n}$, we denote this bimodule by $\mathscr{E}_{n,n-1}$. We can view $\mathscr{E}_{n,n-1}$ as an $(\mathscr{A},\mathscr{A})$-bimodule, so we have an $(\mathscr{A},\mathscr{A})$-bimodule
$$\mathscr{E}:=\bigoplus_{n\geq0}\F_\ell G_ne_{V_nD_n}.$$
Clearly the bimodule $\mathscr{E}$ corresponds to the functor $E$. Similarly, we have an $(\mathscr{A},\mathscr{A})$-bimodule
$$\mathscr{F}:=\bigoplus_{n\geq0}e_{V_nD_n}\F_\ell G_n,$$
which corresponds to the functor $F$. It is obvious that $\mathscr{A}'= k'\otimes_{\F_\ell} \mathscr{A}$ as $k'$-algebras; $\mathscr{E}'= k'\otimes_{\F_\ell}\mathscr{E}$ and $\mathscr{F}'= k'\otimes_{\F_\ell}\mathscr{F}$ as $(\mathscr{A}',\mathscr{A}')$-bimodules.

We identify $\mathscr{A}$ with the subalgebra $1\otimes \mathscr{A}$ of $\mathscr{A}'$ and identify $\mathscr{E}$ with the $(\mathscr{A},\mathscr{A})$-submodule $1\otimes \mathscr{E}$ of $\mathscr{E}'$. By definition, we see that the endomorphism $\mathscr{X}$ of $\mathscr{E}'$ restricts an endomorphism of $\mathscr{E}$. By Proposition \ref{eigenvalue}, any eigenvalue $a$ of $\mathscr{X}$ as an endomorphism of the bimodule $\mathscr{E}'$ are contained in $\F_\ell$.
By elementary linear algebra, if we let $\mathscr{E}_a$ be the generalised $a$-eigenspace of $\mathscr{X}$ on $\mathscr{E}$, then we have $\mathscr{E}'_a= k'\otimes_{\F_\ell}\mathscr{E}_a$. Since $E$ and $F$ are left and right adjoint, there is a corresponding decomposition $F=\bigoplus_{a\in \F_\ell^\times}F_a$, such that $F_a$ is left and right adjoint to $E_a$. Since the bimodule $\mathscr{F}$ corresponds to the functor $F$, there is a decomposition $\mathscr{F}=\bigoplus_{a\in \F_\ell^\times}\mathscr{F}_a$, such that the $(\mathscr{A},\mathscr{A})$-bimodule $\mathscr{F}_a$ corresponds to the functor $F_a$. Since $\mathscr{E}'_a= k'\otimes_{\F_\ell}\mathscr{E}_a$, by the uniqueness of right adjoints (see e.g. \cite[Theorem 2.3.7]{Linckelmann}), we see that $\mathscr{F}'_a= k'\otimes_{\F_\ell}\mathscr{F}_a$.

In \cite[\S6.1]{CR}, the construction of $\Theta'_\lambda$ only uses the functors $E'_{a_j}$, $F'_{a_j}$, the co-unit $\varepsilon'$ of the adjoint pair $(E'_{a_j},F'_{a_j})$ and some elements of the form $c_i^\tau$ in the affine Hecke algebras, where $i$ is some integer and $\tau\in \{1,{\rm sgn}\}$.  (See \cite[\S3.1.1]{CR} for the definition of affine Hecke algebras and see \cite[\S3.1.4]{CR} for the definition of $c_i^\tau$.) We note that in the definition of affine Hecke algebras, the base field can be any field. We also note that the element $c_i^\tau$ is an $\F_\ell$-linear combination of the generators of an affine Hecke algebra, hence $c_i^\tau$ can still be defined even if the base field is $\F_\ell$.

Next, we explain that the complex of functors
$\Theta'_\lambda: {\rm Comp}(\mathcal{A}'_{-\lambda})\to {\rm Comp}(\mathcal{A}'_\lambda)$ in \cite[\S6.1]{CR} can still be defined even if the base field is $\F_\ell$ (in our case). Recall that $(\Theta'_\lambda)^{-r}$ is the restriction of $E_{a_j}'^{({\rm sgn},\lambda+r)}F_{a_j}'^{(1,r)}$ to $\mathcal{A}'_{-\lambda}$ for $r,\lambda+r\geq 0$ and $(\Theta'_\lambda)^{-r}=0$ otherwise. Note that in our case, $\mathcal{A}'_{-\lambda}$ is exactly the category $B'_{j-1}$-mod and $\mathcal{A'}_\lambda$
is exactly the category $B'_j$-mod.  Set $\mathcal{A}_{-\lambda}:=B_{j-1}$-mod and $\mathcal{A}_{\lambda}:=B_j$-mod.  As in \cite[\S6.1]{CR}, we denote by $(\Theta_\lambda)^{-r}$ the restriction of $E_{a_j}^{({\rm sgn},\lambda+r)}F_{a_j}^{(1,r)}$ to $\mathcal{A}_{-\lambda}$ for $r,\lambda+r\geq 0$ and put $(\Theta_\lambda)^{-r}=0$ otherwise.  Since $E_{a_j}'^{({\rm sgn},\lambda+r)}F_{a_j}'^{(1,r)}$ restricts a functor from $\mathcal{A}'_{-\lambda}$ to $\mathcal{A}'_\lambda$, and since we have $\mathscr{E}_{a_j}'= k'\otimes_{\F_\ell}\mathscr{E}_{a_j}$ and $\mathscr{F}_{a_j}'= k'\otimes_{\F_\ell}\mathscr{F}_{a_j}$, $E_{a_j}^{({\rm sgn},\lambda+r)}F_{a_j}^{(1,r)}$ must send an object of $\mathcal{A}_{-\lambda}$ to $\mathcal{A}_\lambda$. In other words, $E_{a_j}^{({\rm sgn},\lambda+r)}F_{a_j}^{(1,r)}$ restricts a functor from $\mathcal{A}_{-\lambda}$ to $\mathcal{A}_{\lambda}$. So $(\Theta_\lambda)^{-r}$ is actually a functor from $\mathcal{A}_{-\lambda}$ to $\mathcal{A}_\lambda$.

In the third paragraph of \cite[\S6.1]{CR}, Chuang and Rouquier defined a map $d'^{-r}:(\Theta'_\lambda)^{-r}\to (\Theta'_\lambda)^{-r+1}$. Using the same way (replacing the co-unit $\varepsilon'$ of the adjoint pair $(E'_{a_j},F'_{a_j})$ by the co-unit $\varepsilon$ of the adjoint pair $(E_{a_j},F_{a_j})$), we can define a map $d^{-r}:(\Theta_\lambda)^{-r}\to (\Theta_\lambda)^{-r+1}$.
So we obtain a complex of functors $$\Theta_\lambda=\cdots\to(\Theta_\lambda)^i\xrightarrow{d^i}(\Theta_\lambda)^{i+1}\to\cdots.$$
Evaluating $\Theta_\lambda$ at the $(B_{j-1},B_{j-1})$-bimodule $B_{j-1}$, we obtain a complex
$$C_j:=\cdots\to(\Theta_\lambda)^i(B_{j-1}) \xrightarrow{d^i(B_{j-1})}(\Theta_\lambda)^{i+1}(B_{j-1})\to\cdots.$$
Since we have $\mathscr{E}_{a_j}'= k'\otimes_{\F_\ell}\mathscr{E}_{a_j}$ and $\mathscr{F}_{a_j}'= k'\otimes_{\F_\ell}\mathscr{F}_{a_j}$ and since the co-unit of the adjoint pair $(E'_{n,n-1},F'_{n,n-1})$ and the co-unit of the adjoint pair $(E_{n,n-1},F_{n,n-1})$ are constructed in the same way (i.e., the construction does not depend on the base field), it is easy to see that
$C'_j\cong k'\otimes_{\F_\ell}C_j$ as complexes of $(B'_j,B'_{j-1})$-bimodules. $\hfill\square$

\section{Proof of Theorem \ref{main}}\label{s4}

In Section \ref{s3}, we showed that the splendid Rickard equivalence in Remark \ref{1.4} (iii) descends to $\F_\ell$, so we finished half of the task listed in Remark \ref{1.5}. To finish the proof of Theorem \ref{main}, we need to show that the splendid Rickard equivalences in Remark \ref{1.4} (i) and (ii) descend to $\F_\ell$.

Keep the notation of Remark \ref{1.4} and \ref{1.5}. Recall that $G$ is $\GL_n(q)$, $b$ is a unipotent block of $k'G$, $P$ is a defect group of $b$, $c$ is the block of $k'N_G(P)$ corresponding to $b$ via the Brauer correspondence, $e$ is the multiplicative order of $q$ in $k^\times$, and $w$ is the ($e$-)weight of $b$. The block $b$ corresponds uniquely to an $e$-core $\kappa$. The assumption that the defect group $P$ of $b$ is abelian forces $w<\ell$. Let $m$ be the greatest integer such that $\ell^m$ divides $q^e-1$, then $P\cong \underbrace{C_{\ell^m}\times\cdots\times C_{\ell^m}}_w$ (see \cite[\S1.5]{Turner}). (Here the notation $C_{\ell^m}$ denotes a cyclic group of order $\ell^m$.) We identify $P$ and $\underbrace{C_{\ell^m}\times\cdots\times C_{\ell^m}}_w$ via the isomorphism. Let $t:=n-ew$. Then by \cite[\S 1.4 and \S1.5]{Turner}, the block algebra $k'N_G(P)c$ is isomorphic to $k' N_{\GL_{ew}(q)}(P)c_{ew}\otimes_{k'} k' \GL_t(q)c_0$, where $c_{ew}$ is the principal block of $k'N_{\GL_{ew}(q)}$ and $c_0$ is the unipotent block of $k'\GL_t(q)$ corresponding to the $e$-core $\kappa$ and having defect group $1$. We identify $k'N_G(P)c$ and $k' N_{\GL_{ew}(q)}(P)\otimes_{k'} k' \GL_t(q)c_0$ via the isomorphism. Since $c\in \F_\ell N_G(P)$, $c_{ew}\in\F_\ell N_{\GL_{ew}(q)}(P)$ and $c_0\in \F_\ell\GL_t(q)$ (see Corollary \ref{unipotent} and \ref{principal}), we also have $\F_\ell N_G(P)c\cong \F_\ell N_{\GL_{ew}(q)}(P)c_{ew}\otimes_{\F_\ell} \F_\ell \GL_t(q)c_0$. We identify these two algebras.

By \cite[Lemma 6]{CK}, there is a splendid Morita equivalence between $k'N_{\GL_{ew}(q)}(P)c_{ew}$ and $k'N_G(P)c$. To show that this splendid Morita equivalence descends to $\F_\ell$, we prove a descent proposition for \cite[Lemma 6]{CK}.

\begin{proposition}\label{defectzero}
Let $G_1$ and $G_2$ be finite groups. Let $b_1$ and $b_2$ be blocks of $k'G_1$ and $k'G_2$, and assume that $b_2$ has defect group $1$. Assume that $k$ is a subfield of $k'$ such that $b_1\in kG_1$ and $b_2\in k G_2$, then $kG_1b_1$ and $kG_1b_1\otimes_k kG_2b_2$ (a block algebra of $G_1\times G_2$) are splendidly Morita equivalent.
\end{proposition}

\noindent{\it Proof.} Since every finite group has a finite splitting field, we may assume that $k'$ is finite. Let $i$ be a primitive idempotent in $k'G_2b_2$. By the proof of \cite[Lemma 6]{CK}, the $(k'(G_1\times G_2), k'G_1)$-bimodule $k'G_1b_1\otimes_{k'} k'G_2i$ induces a splendid Morita equivalence between $k'G_1b_1\otimes_{k'}k'G_2b_2$ and $k'G_1b_1$. Let $\Gamma:={\rm Gal}(k'/k)$. Since $b_2$ has defect group $1$, $k'G_2i$ is the unique projective $k'Gb_2$-module, up to isomorphism. For any $\sigma\in \Gamma$, ${}^\sigma (k'G_2i)$ is still a projective $k'G_2b_2$-module, so we have ${}^\sigma (k'G_2i)\cong k'G_2i$. Then by \cite[Lemma 6.2 (c)]{Kessar_Linckelmann}, there is a projective $kG_2b_2$-module $Y$ such that $k'G_2i\cong k'\otimes_k Y$. It follows that
$$k'G_1b_1\otimes_{k'} k'G_2i\cong k'\otimes_k (kG_1b_1\otimes_{k} Y)$$
as $(k'(G_1\times G_2)(b_1\otimes b_2), k'G_1b_1)$-bimodule. By \cite[Proposition 4.5 (c)]{Kessar_Linckelmann}, the $(k(G_1\times G_2)(b_1\otimes b_2), kG_1b_1)$-bimodule $kG_1b_1\otimes_{k} Y$ induces a Morita equivalence between $kG_1b_1\otimes_kkG_2b_2$ and $kG_1b_1$. By \cite[Lemma 5.1 and 5.2]{Kessar_Linckelmann}, this Morita equivalence is splendid.  $\hfill\square$

\medskip By Proposition \ref{defectzero}, $\F_\ell N_{\GL_{ew}(q)}(P)c_{ew}$ and $\F_\ell N_G(P)c$ are splendidly Morita equivalent. By \cite[\S1.4]{Turner}, we have $N_{\GL_{ew}(q)}(P)\cong N_{\GL_e(q)}(C_{\ell^m})\wr S_w$.
Let $b_e$ be the principal block of $k'\GL_e(q)$, then $b_e$ is the unique unipotent block of $k'\GL_e(q)$ corresponding to the empty $e$-core. So $b_e$ has $C_{\ell^m}$ as a defect group. The Brauer correspondent of $b_e$ in $k'N_{\GL_e(q)}(C_{\ell^m})$ is the principal block of $k'N_{\GL_e(q)}(C_{\ell^m})$, and we denote it by $c_e$. By Rouquier's result on cyclic blocks (\cite[Theorem 10.1]{Rou}), there is a complex of $C'$ of $(k'\GL_e(q)b_e,k'N_{\GL_e(q)}(C_{\ell^m})c_e)$-bimodules inducing a splendid Rickard equivalence between $k'\GL_e(q)b_e$ and $k'N_{\GL_e(q)}(C_{\ell^m})c_e$. Then by a theorem of Marcus (\cite[Theorem 4.3]{Marcusonequivalences}), the complex $C'\wr S_w$ induces a splendid Rickard equivalence between the principal block algebra of $k'(\GL_e(q)\wr S_w)$ and the principal block algebra of $k'(N_{\GL_e(q)}(C_{\ell^m})\wr S_w)$. By \cite[Theorem 1.10]{Kessar_Linckelmann}, Rouquier's complex $C'$ descends to $\F_\ell$, namely that there is a Rickard complex $C$ of $(\F_\ell\GL_e(q)b_e,\F_\ell N_{\GL_e(q)}(C_{\ell^m})c_e)$-bimodules, such that $C'\cong k'\otimes_{\F_\ell}C$. We note that the blanket assumption in [19] that the coefficient field is ``big enough" is not used in [19, Theorem 4.3] and its proof. Hence by \cite[Theorem 4.3]{Marcusonequivalences}, the complex $C\wr S_w$ induces a splendid Rickard equivalence between the principal block algebra of $\F_\ell(\GL_e(q)\wr S_w)$ and the principal block algebra of $\F_\ell(N_{\GL_e(q)}(C_{\ell^m})\wr S_w)$. Now we have proved the following statement.

\begin{proposition}\label{4.2}
$\F_\ell N_G(P)c$ is splendidly Morita equivalent to the principal block algebra of $\F_\ell(\GL_e(q)\wr S_w)$.
\end{proposition}

Consider an abacus having $w+i(w-1)$ beads on the $i$-th runner, where $i\in\{0,1,\cdots,e-1\}$. Let $\rho$ be the $e$-core have this abacus representation, let $r:=|\rho|$, and let $n':=we+r$. By results of Miyachi (\cite[Theorem 5.0.7]{Mi01}) and Turner (\cite[Theorem 1]{Turner}), the unipotent block (say $b'$) of $k'\GL_{n'}(q)$ corresponding to the $e$-core $\rho$ is splendidly Morita equivalent to the principal block algebra of $k'(\GL_e(q)\wr S_w)$. The integer $n'$ is exactly the integer $n'$ mentioned in Remark \ref{1.4}, and the block $b'$ of $k'\GL_{n'}(q)$ is exactly the unipotent block $b'$ mentioned in Remark \ref{1.4}.

Let us review the construction of the splendid Morita equivalence constructed by Turner. Let $f_0$ be the unipotent block of $k'\GL_r(q)$ corresponding to the $e$-core $\rho$, so $f_0$ has defect group $1$.
Denote by $b_*$ the principal block of $k'(\GL_e(q)\wr S_w)$. By \cite[Lemma 6]{CK}, $k'(\GL_e(q)\wr S_w)b_*$ is splendidly Morita equivalent to $k'(\GL_e(q)\wr S_w)b_*\otimes_{k'} k'\GL_r(q)f_0$. Let $L$ be the Levi subgroup
$\underbrace{\GL_e(q)\times \cdots \times \GL_e(q)}_w\times \GL_r(q)$
of $\GL_{n'}(q)$. Recall that $b_e$ denotes the principal block of $k'\GL_e(q)$. Let
$f:=\underbrace{b_e\otimes\cdots\otimes b_e}_w\otimes f_0$, a block idempotent of $k'L$. We can view $S_w$ as the subgroup of permutation matrices of $\GL_{n'}(q)$ whose conjugation action on $L$ permutes the factors of $\underbrace{\GL_e(q)\times \cdots \times \GL_e(q)}_w$. Let $N$ be the semi-direct product of $L$ and $S_w$, a subgroup of $\GL_{n'}(q)$ isomorphic to $\GL_e(q)\wr S_w\times\GL_r(q)$. Clearly the conjugation action of $N$ stabilises $f$, hence $f$ is an idempotent in the center of $k'N$. By \cite[Lemma 1 (3)]{Turner}, $f$ is also a block of $k'N$. This forces that the idempotent $\underbrace{b_e\otimes\cdots\otimes b_e}_w$ is a block of $k'(\GL_e(q)\wr S_w)$, and hence $\underbrace{b_e\otimes\cdots\otimes b_e}_w$ must equal to the principal block $b_*$ of $k'(\GL_e(q)\wr S_w)$. So the block algebra $k'Nf$ is isomorphic to the algebra  $k'(\GL_e(q)\wr S_w)b_*\otimes_{k'} k'\GL_r(q)f_0$. We identify these two algebras via the isomorphism. Let $D:=\underbrace{C_{\ell^m}\times\cdots\times C_{\ell^m}}_w$, a Sylow $\ell$-subgroup of $\underbrace{\GL_e(q)\times \cdots \times \GL_e(q)}_w$. By \cite[Lemma 1 (4)]{Turner}, $k'\GL_{n'}(q)b'$ and $k'Nf$ both have defect group $D$ and are Brauer correspondents.

By Alperin's description of the Brauer correspondence (\cite[Lemma 6.2.7]{Benson}), the $k'(\GL_{n'}(q)\times \GL_{n'}(q)^{\rm op})$-module $k'\GL_{n'}(q)b'$ and the $k'(N\times N^{\rm op})$-module $k'Nf$ both have vertex $\Delta D$ and are Green correspondents. Let $T'$ be the Green correspondent of $k'\GL_{n'}(q)b'$ in $\GL_{n'}(q)\times N^{\rm op}$, an indecomposable summand of ${\rm Res}_{\GL_{n'}(q)\times N^{\rm op}}^{\GL_{n'}(q)\times \GL_{n'}(q)^{\rm op}}(k'\GL_{n'}(q)b')$. Since $k'Nf$ is a direct summand of ${\rm Res}_{N\times N^{\rm op}}^{\GL_{n'}(q)\times N^{\rm op}}(T')$, we have $T'f\neq 0$, thus $T'f=T$ and $T'$ is a $(k'\GL_{n'}(q)b',k'Nf)$-bimodule. By \cite[Proposition 1]{Turner}, the $(k'\GL_{n'}(q)b',k'Nf)$-bimodule $T'$ induces a splendid Morita equivalence between $k'\GL_{n'}(q)b'$ and $k'Nf$. So there is a splendid Morita equivalence between $k'\GL_{n'}(q)b'$ and $k'(\GL_e(q)\wr S_w)b_*$. Next, we prove that this splendid Morita equivalence descends to $\F_\ell$.

\begin{proposition}\label{4.3}
There is a splendid Morita equivalence between $\F_\ell\GL_{n'}(q)b'$ and the principal block algebra of $\F_\ell(\GL_e(q)\wr S_w)$.
\end{proposition}

\noindent{\it Proof.} By Proposition \ref{defectzero}, $\F_\ell(\GL_e(q)\wr S_w)b_*$ is splendidly Morita equivalent to $\F_\ell(\GL_e(q)\wr S_w)b_*\otimes_{\F_\ell} \F_\ell\GL_r(q)f_0$.  Since $b_e\in \F_\ell\GL_e(q)$ and $f_0\in \F_\ell \GL_r(q)$ (see Corollary \ref{unipotent}), we have $f=\underbrace{b_e\otimes\cdots\otimes b_e}_w\otimes f_0\in \F_\ell N$. Since the block algebra $k'Nf$ is isomorphic to $k'(\GL_e(q)\wr S_w)b_*\otimes_{k'} k'\GL_r(q)f_0$, the block algebra $\F_\ell Nf$ is isomorphic to $\F_\ell(\GL_e(q)\wr S_w)b_*\otimes_{k'} \F_\ell\GL_r(q)f_0$. Hence $\F_\ell(\GL_e(q)\wr S_w)b_*$ is splendid Morita equivalent to $\F_\ell Nf$.

By the definition of Green correspondent, $T'$ is the unique (up to isomorphism) direct summand of the $(k'\GL_{n'}(q)b',k'Nf)$-bimodule $k'\GL_{n'}(q)b'$ having $\Delta D$ as a vertex. Noting that $k'\GL_{n'}(q)b'= k'\otimes_{\F_\ell} \F_\ell \GL_{n'}(q)b'$, by \cite[Lemma 5.1]{Kessar_Linckelmann}, there is an indecomposable direct summand $T$ of the $(\F_\ell \GL_{n'}(q)b',\F_\ell Nf)$-bimodule $\F_\ell\GL_{n'}(q)b'$ such that $T'\cong k'\otimes_{\F_\ell} T$ (by the uniqueness of $T'$). By \cite[Proposition 4.5]{Kessar_Linckelmann}, $T$ induces a splendid Morita equivalence between $\F_\ell \GL_{n'}(q)b'$ and $\F_\ell Nf$. This completes the proof. $\hfill\square$

\medskip\noindent{\it Proof of Theorem \ref{main}.} By Proposition \ref{3.3}, \ref{4.2} and \ref{4.3}, the task listed in Remark \ref{1.5} is finished, hence the proof of Theorem \ref{main} is complete. $\hfill\square$

\section{On general blocks of $\GL_n(q)$}\label{ap}

By the proof of \cite[Theorem 1.12]{Kessar_Linckelmann}, to answer the question whether Conjecture \ref{KLconj} holds for general blocks of $\GL_n(q)$, it suffices to answer the following question.

\begin{question}\label{A1}
Let $G$ be the group $\GL_n(q)$, and let $b$ be a block of $\O'G$ with an abelian defect group $P$. Let $c$ be the block of $\O'N_G(P)$ corresponding to $b$ via the Brauer correspondence. Suppose that $b\in \O G$. Then $c\in \O N_G(P)$. Are the block algebras $\O Gb$ and $\O N_G(P)c$ splendidly Rickard equivalent? 
\end{question}

By the proof of \cite[Theorem 7.20]{CR}, Chuang and Rouquier first reduced the statement in \cite[Theorem 7.20]{CR} to unipotent blocks by using results in \cite{BDR}. More precisely, by results in \cite{BDR}, there is a complex $X'_1$ inducing a splendid Rickard equivalence between $\O'Gb$ and a unipotent block (say $b_1$) of $\O'G_1$, where $G_1\cong\GL_{n_1}(q^{d_1})\times\cdots\times\GL_{n_r}(q^{d_r})$ is a subgroup of $G$. Let $D$ be a defect group of $b_1$. By \cite[Theorem 1.15]{Ruh}, $D$ is also a defect group of $b$. We may assume that $P=D$ (because we can change the choice of $P$).  Let $c_1$ be the block of $\O'N_{G_1}(P)$ corresponding to $b_1$ via the Brauer correspondence. Then the complex
$$Y'_1:={\rm Ind}_{N_{G\times G_1^{\rm op}}(\Delta P)}^{N_G(P)\times N_{G_1}(P)^{\rm op}}{\rm Br}_{\Delta P}(X'_1)$$
(where ${\rm Br}_{\Delta P}(X')$ denotes the $\Delta P$-Brauer construction of $X'$) induces a splendid Rickard equivalence between $\O'N_G(P)c$ and $\O'N_{G_1}(P)c_1$ (see \cite[Proposition 1.36 and Remark 1.37]{Ruh}). Hence the statement in \cite[Theorem 7.20]{CR} reduces to proving that $\O'G_1b_1$ and $\O'N_{G_1}(P)c_1$ are splendidly Rickard equivalent. In other words, the statement in \cite[Theorem 7.20]{CR} is reduced to unipotent blocks.

By Corollary \ref{unipotent}, we have $b_1\in \Z_\ell G_1\subseteq \O G_1$, and hence we have $c_1\in \Z_\ell N_{G_1}(P)\subseteq  \O N_{G_1}(P)$. Since we have answered Question \ref{A1} for unipotent blocks of general linear groups (Theorem \ref{main}), if we can prove that there is a complex $X_1$ of $(\O Gb,\O G_1b_1)$-bimodules such that $X'_1\cong \O'\otimes_{\O}X_1$, then Question \ref{A1} has a positive answer. We will give a sufficient condition under which there exists such a complex $X_1$ (see Theorem \ref{sufficient} below). Then by \cite[Proposition 4.5 (a)]{Kessar_Linckelmann}, $X_1$ induces a splendid Rickard equivalence between $\O Gb$ and $\O G_1b_1$, and $Y_1$ induces a splendid Rickard equivalence between $\O N_G(P)c$ and $\O N_{G_1}(P)c_1$, where
$$Y_1:={\rm Ind}_{N_{G\times G_1^{\rm op}}(\Delta P)}^{N_G(P)\times N_{G_1}(P)^{\rm op}}{\rm Br}_{\Delta P}(X_1).$$

Let us review the construction of the Rickard complex $X'_1$ constructed by Bonnaf\'{e}, Dat and Rouquier \cite{BDR}. So we should first review some material in \cite{BDR}.

\subsection{The Deligne--Lusztig induction}

Assume that $\Lambda\in \{\O,k\}$. Let $\mathbf{G}$ be a connected reductive algebraic group over an algebraic closure of a finite field whose characteristic is not $\ell$, endowed with a Frobenius endomorphism $F$. Let $\L$ be an $F$-stable Levi subgroup of $\G$ contained in a parabolic subgroup $\mathbf{P}$ with unipotent radical $\mathbf{V}$ such that $\mathbf{P}=\V\rtimes \L$. The Deligne--Lusztig variety
$$\Y_\P:=\{g\V\in \G / \V~|~g^{-1}F(g)\in \V\cdot F(\V)\}$$
has a left action of $\G^F$ and a right action of $\L^F$ by multiplication. By works of Rickard (\cite{Rickard94}) and Rouquier (\cite{Rou2002etalas}), there is an object $\GG_c(\Y_\P,\Lambda)$ of ${\rm Ho}^b(\Lambda (\G^F\times (\L^F)^{\rm op})\mbox{-perm})$ associated with $\mathbf{Y}_\P$, well defined up to isomorphism, where $\Lambda (\G^F\times (\L^F)^{\rm op})\mbox{-perm}$ denotes the category of finitely generated $\ell$-permutation $(\Lambda\G^F,\Lambda\L^F)$-bimodules.

\begin{proposition}\label{etaladescent}
$\GG_c(\Y_\P,k)\cong k\otimes_{\F_\ell} \GG_c(\Y_\P,\F_\ell)$ as complexes of $(k\G^F,k\L^F)$-bimodules.
\end{proposition}

\noindent{\it Proof.} The statement follows from \cite[Lemma 2.8]{Rickard94}: take the module ``$U$" in \cite[Lemma 2.8]{Rickard94} to be the direct sum of all the permutation $\F_\ell(\G^F\times (\L^F)^{\rm op})$-modules of the form $\F_\ell(\G^F\times (\L^F)^{\rm op} /H)$, where $H$ is a subgroup of $\G^F\times (\L^F)^{\rm op}$. Then the category of finitely generated $\ell$-permutation $\F_\ell(\G^F\times (\L^F)^{\rm op})$-modules is exactly the category ``${\rm add}$-$U$" in \cite[Lemma 2.8]{Rickard94}. $\hfill\square$

\medskip Let $\Lambda\in \{K',k'\}$. The complex $\RG_c(\Y_\P,\Lambda)$ of $\ell$-adic cohomology induces a morphism
$$R_{\L\subset \P}^\G:G_0(\Lambda\L^F)\to G_0(\Lambda\G^F)$$
between Grothendieck groups, which is called the {\it Deligne--Lusztig induction}.

\subsection{The Jordan decomposition}\label{subsection:The Jordan decomposition}

Let $\mathbf{G}^*$ be a group dual to $\mathbf{G}$ with Frobenius endomorphism $F^*$.
Let $\Irr_{K'}(\mathbf{G}^F)$ denote the set of characters of irreducible representations of $\mathbf{G}^F$ over $K'$. Deligne and Lusztig gave a decomposition of $\Irr_{K'}(\mathbf{G}^F)$ into rational series
$$\Irr_{K'}(\mathbf{G}^F)=\coprod_{(s)}\Irr_{K'}(\mathbf{G}^F,(s)),$$
where $(s)$ runs over the set of $\mathbf{G}^*{}^{F^*}$-conjugacy classes of semi-simple elements of $\mathbf{G}^*{}^{F^*}$. The unipotent characters of $\mathbf{G}^F$ are those in $\Irr_{K'}(\mathbf{G}^F,(1))$.

Let $s$ be a semi-simple element of ${\G^*}^{F^*}$ of order prime to $\ell$.  Brou\'{e} and Michel (\cite{BM}) showed that $\coprod_{(t)}\Irr_{K'}(\mathbf{G}^F,(t))$, where $(t)$ runs over conjugacy classes of semi-simple elements of $\mathbf{G}^*{}^{F^*}$ whose $\ell'$-part is $(s)$, is a union of blocks of $\O'\mathbf{G}^F$. The sum of the corresponding block idempotents is an idempotent $e_s^{\mathbf{G}^F}$ in the center of $\O'\mathbf{G}^F$, so there is a decomposition $1=\sum_{(s)}e_s^{\G^F}$, where $(s)$ runs over $\mathbf{G}^*{}^{F^*}$-conjugacy classes of semi-simple $\ell'$-elements of $\mathbf{G}^*{}^{F^*}$. We denote the set $\coprod_{(t)}\Irr_{K'}(\mathbf{G}^F,(t))$ above by ${\rm Irr}_{K'}(\G^F,e_s^{\G^F})$.

Let $\L$ be an $F$-stable Levi subgroup of $\G$ with dual $\L^*\subset \G^*$ containing $C_{\G^*}(s)$. By a result of Lusztig (see \cite[Theorem 11.4.3]{DM}), there is a sign $\varepsilon_{\L,\G}\in\{1,-1\}$ such that $\varepsilon_{\L,\G} R_{\L\subset \P}^{\G}$ induces a bijection
$$\Irr_{K'}(\L^F,(s))\xrightarrow{\sim}\Irr_{K'}(\G^F,(s)),~\psi\mapsto \varepsilon_{\L,\G} R_{\L\subset \P}^\G(\psi).$$
In this situation with $s$ being an $\ell'$-element, $\varepsilon_{\L,\G} R_{\L\subset \P}^{\G}$ also induces a bijection
$$\Irr_{K'}(\L^F,e_s^{\L^F})\xrightarrow{\sim}\Irr_{K'}(\G^F,e_s^{\G^F}),~\psi\mapsto \varepsilon_{\L,\G} R_{\L\subset \P}^\G(\psi).$$

\subsection{A sufficient condition}

Denote by $\bar{\F}_q$ the algebraic closure of $\F_q$.
Let $\mathbf{G}:=\GL_n(\bar{\F}_q)$, and let $F$ be the map $\mathbf{G}\to \mathbf{G}$ sending every $A=(a_{ij})_{1\leq i,j\leq n}\in \GL_n(\bar{\F}_q)$ to $(a_{ij}^q)_{1\leq i,j\leq n}$. Hence $\mathbf{G}^F$ is $G:=\GL_n(q)$. The pair $(\mathbf{G},F)$ is dual to itself (see e.g. \cite[Examples 11.1.13]{DM}).

We return to the context of Question \ref{A1} and keep the notation of the beginning part of this section. We will express a complex $X_1'$ more explicitly. Since $b$ is a block of $\O'G=\O'\GL_n(q)=\O'\G^F$, there is a semi-simple $\ell'$-element $s$ of $\G^F$ such that $be_s^{\G^F}=b=e_s^{\G^F}b$. Let $\mathbf{L}:=C_{\G}(s)$, by \cite[Theorem 11.7.3]{DM}, $\L$ is a Levi subgroup of a parabolic subgroup $\mathbf{P}$ of $\G$. Let $\mathbf{V}$ be the unipotent radical of $\mathbf{P}$, then we have $\mathbf{P}=\mathbf{V}\rtimes \mathbf{L}$.
Clearly $\L$ is $F$-stable. The group $G_1$ is exactly the group $\L^F$.

Denote by $C'$ the complex $\GG_c(\Y_\P,\O')$ of $(\O'G,\O'G_1)$-bimodules. By \cite[Theorem 7.7]{BDR}, there is a block $b'$ of $\O'G_1$ such that the complex $b{C'}^{\rm red}b'$ induces a splendid Rickard equivalence between $\O'Gb$ and $\O'G_1b'$. (Note that in our case, $\L^F$ and $\mathbf{N}^F$ in \cite[Theorem 7.7]{BDR} are equal.) It is clear that $bC'b'$ is isomorphic to $bC'^{\rm red}b'$ plus a complex of $(\O' Gb,\O'G_1b')$-bimodules which is homotopy equivalent to zero. Then by the definition of a splendid Rickard complex, we see that the complex $bC'b'$ as well induces a splendid Rickard equivalence between $\O'Gb$ and $\O'G_1b'$. 

Since $s\in Z(\L)$, there is a bijection
$$\Irr_{K'}(\L^F,e_1^{\L^F})\xrightarrow{\sim}\Irr_{K'}(\L^F,e_s^{\L^F}), \psi\to \eta\psi$$
where $\eta$ is the one-dimensional character of $\L^F$ corresponding to $s$ (see e.g. \cite[Proposition 8.26]{CE}). Let $S'$ an $\O'\L^F$-module which affords the character $\eta$. Then it is easy to see that
the functor (say $\Phi$) sending an $\O'\L^Fe_1^{\L^F}$-module $V$ to the $\O'\L^Fe_s^{\L^F}$-module $V\otimes_{\O'}S'$ induces a Morita equivalence between $\O'\L^Fe_s^{\L^F}$ and $\O'\L^Fe_1^{\L^F}$. Assume that this Morita equivalence is induced by an $(\O'\L^Fe_s^{\L^F},\O'\L^Fe_1^{\L^F})$-bimodule $M'$. Then as  $(\O'\L^Fe_s^{\L^F},\O'\L^Fe_1^{\L^F})$-bimodules,
$$M'\cong M'\otimes_{\O'\L^Fe_1^{\L^F}}\O'\L^Fe_1^{\L^F}\cong \Phi(\O'\L^Fe_1^{\L^F})= \O'\L^Fe_1^{\L^F}\otimes_{\O'} S'.$$
Since $\O'\L^Fe_s^{\L^F}$ and $\O'\L^Fe_1^{\L^F}$ are direct sums of block algebras,  there is a unique block $b_1$ of $\O'\L^F=\O'G_1$ such that $b'M'\cong M'b_1$, and $b'M'b_1$ induces a Morita equivalence between $\O'G_1b'$ and $\O'G_1b_1$. (The $b_1$ in the beginning part of this section is the $b_1$ here.) So the complex
$bC'b'\otimes_{\O'G_1b'} b'M'b_1= bC'\otimes_{\O'G_1} M'b_1$ induces a Rickard equivalence between $\O'Gb$ and $\O'G_1b_1$. Set $\bar{C}':=\GG_c(\Y_\P,k')\cong k'\otimes_{\O'}C'$, $\bar{S}':=k'\otimes_{\O'}S$ and $\bar{M}':=k'\otimes_{\O'}M'$.  Denote by $\bar{b}$ (resp. $\bar{b_1}$) the image of $b$ (resp. $b_1$) in $k'G$ (resp. $k'G_1$).

\begin{proposition}\label{x1}
The complex $\bar{b}\bar{C}'\otimes_{k'G_1} \bar{M}'\bar{b}_1$ is a splendid Rickard complex for $k'G\bar{b}$ and $k'G_1\bar{b}_1$, and it lifts uniquely (up to isomorphism) to a splendid Rickard complex $X'_1$ for $\O'Gb$ and $\O'G_1b_1$.
\end{proposition}

\noindent{\it Proof.} Since $M'\cong \O'\L^F e_1^{\L^F}\otimes_{\O'}S'$, $M'$ is isomorphic to a direct summand of
$$\O'\L^{F}\otimes_{\O'}S'=\O'G_1\otimes_{\O'}S'\cong {\rm Ind}_{\Delta G_1}^{G_1\times G_1^{\rm op}}(S'),$$
where the isomorphism is by \cite[Proposition 2.8.19]{Linckelmann}.
It follows that $\bar{M}'$ is isomorphic to a direct summand of ${\rm Ind}_{\Delta G_1}^{G_1\times G_1^{\rm op}}(\bar{S}')$. Note that $\bar{S}'$ is a 1-dimension $k'G_1$-module, hence the restriction of $\bar{S}'$ to any $p$-subgroup of $\Delta G_1$ is a trivial module. Using the Mackey formula, we easily see that $\bar{M}'$ is an $\ell$-permutation $k'(G_1\times G_1^{\rm op})$-module. So the complex $\bar{b}\bar{C}'\otimes_{k'G_1} \bar{M}'\bar{b}_1$ is a complex of $\ell$-permutation $k'(G\times G_1^{\rm op})$-modules, hence it induces a splendid Rickard equivalence between $k'G\bar{b}$ and $k'G_1\bar{b}_1$. By Theorem \ref{lifting}, there is a unique (up to isomorphism) complex $X'_1$ of $(\O'Gb,\O'G_1b_1)$-bimodules inducing a splendid Rickard equivalence between $\O'Gb$ and $\O'G_1b_1$ such that $k'\otimes_{\O'} X'_1\cong \bar{b}\bar{C}'\otimes_{k'G_1} \bar{M}'\bar{b}_1$. $\hfill\square$

\medskip In Question \ref{A1}, we assumed that $b\in \O G$. Under this assumption, we don't know whether the values of the character $\eta$ are contained in $\O$.  If we assume that the values of $\eta$ are contained in $\O$, then we can give a positive answer to Question \ref{A1} for the block $b$.

\begin{theorem}\label{sufficient}
Keep the notation above. Assume that the values of  $\eta$ are contained in $\O$, then there is a complex $X'_1$ of $(\O Gb,\O G_1b_1)$-bimodules inducing a Rickard equivalence between $\O Gb$ and $\O G_1b_1$, and satisfying $X'_1\cong \O'\otimes_\O X_1$.
\end{theorem}

\noindent{\it Proof.} Let $\bar{C}:=\GG_c(\Y_\P,k)$, a complex of $(kG,kG_1)$-bimodule. By Proposition \ref{etaladescent}, we have $\bar{C}'\cong k'\otimes_k \bar{C}$. Hence $\bar{b}\bar{C}'\cong k'\otimes_k \bar{b}\bar{C}$ as complexes of $(k'Gb, k'G_1)$-bimodules. Let $\bar{\eta}$ be the $k'$-character afforded by the $kG_1$-module $\bar{S}$. By \cite[Proposition 4.3.5]{Linckelmann}, the values of  $\eta$ are contained in $\O$ if and only if  the values of $\bar{\eta}$ are contained in $k$. So by assumption, there exists a $kG_1$-module $\bar{S}$ such that $\bar{S}'\cong k'\otimes_k \bar{S}$. Since $e_1^{\L^F}$ is a sum of unipotent blocks of $\O'\L^F=\O'G_1$, by Corollary \ref{unipotent}, we have $e_1^{\L^F}\in \Z_\ell G_1 \subseteq\O G_1$. Let $\bar{e}_1^{\L^F}$ be the image of $e_1^{\L^F}$ in $k'G_1$, then we have $\bar{e}_1^{\L^F}\in kG_1$. Taking $\bar{M}:=k\L^F\bar{e}_1^{\L^F}\otimes_k \bar{S}$, we see that $\bar{M}'\cong k'\otimes_k \bar{M}$. So we have
$$\bar{b}\bar{C}'\otimes_{k'G_1} \bar{M}'\bar{b}_1\cong k'\otimes_k (\bar{b}\bar{C}\otimes_{k G_1} \bar{M}\bar{b}_1),$$
as complexes of $(k'G,k'G_1)$-bimodules. The same argument in the proof of Proposition \ref{x1} shows that $\bar{b}\bar{C}\otimes_{k G_1} \bar{M}\bar{b}_1$ is a complex of $\ell$-permutation $k(G\times G_1^{\rm op})$-modules. By \cite[Proposition 4.5 (a)]{Kessar_Linckelmann}, $\bar{b}\bar{C}\otimes_{k G_1} \bar{M}\bar{b}_1$ induces a splendid Rickard equivalence between $kGb$ and $kG_1b_1$. By Theorem \ref{lifting},  there is a unique (up to isomorphism) complex $X_1$ of $(\O G_1b,\O G_1b_1)$-bimodules inducing a splendid Rickard equivalence between $\O Gb$ and $\O G_1b_1$ such that $k\otimes_{\O} X_1\cong \bar{b}\bar{C}\otimes_{kG_1} \bar{M}\bar{b}_1$. Since
$$k'\otimes_{\O'} (\O'\otimes_{\O} X_1)\cong k'\otimes_k (k\otimes_{\O}X_1) \cong \bar{b}\bar{C}'\otimes_{k'G_1} \bar{M}'\bar{b}_1\cong k'\otimes_{\O'}X'_1,$$
by \ref{lifting}, we have $\O'\otimes_{\O} X_1\cong X'_1$.  $\hfill\square$

\medskip For unipotent blocks, the hypothesis of Theorem \ref{sufficient} holds. We end this paper by providing an example where the hypothesis of Theorem \ref{sufficient} does not hold. 

\begin{example}\label{example}
{\rm Let $n=2$, $q=5$ and $\ell=3$. So $G=\G^F=\GL_2(5)$. Let $a$ be a generator of $\F_5^\times$, and let $s=\left({\begin{array}{*{20}{c}}
		a&  \\
		& a^{-1}
\end{array}}\right)$, then $s$ is a semi-simple element of $\G^F={\G^*}^{F^*}$ of order prime to $\ell$. Hence $G_1=\L^F=C_{\G^F}(s)=\left\{\left( {\begin{array}{*{20}{c}}
x& \\
&y
\end{array}} \right) \middle|~x,y\in \F_5^\times\right\}$. Since $s$ is of order $4$,  the one-dimensional character $\eta$ of $\L^F$ corresponding to $s$ is of order $4$ as well. So the values of $\eta$ contain a primitive $4$-th root of unity in $K'$. Note that $\L^F$ is abelian, and any element of $\L^F$ is a semi-simple element of order prime to $\ell$. Hence the set ${\rm Irr}_{K'}(\L^F,e_s^{\L^F})$ contains only one character. It follows that the set ${\rm Irr}_{K'}(\G^F,e_s^{\G^F})$ contains only one character (see \S\ref{subsection:The Jordan decomposition}), and hence $e_s^{\G^F}$ is exactly a block of $G$. We denote this block by $b$. Let $\sigma$ be any field automorphism of $K'$ sending any $\ell'$-th root of unity $\epsilon\in K'$ to $\epsilon^\ell$. Then by \cite[Lemma 3.1]{Ruh21}, ${}^\sigma{\rm Irr}_{K'}(\G^F,e_s^{\G^F})={\rm Irr}_{K'}(\G^F,e_{s^\ell}^{\G^F})$. Since $s$ and $s^\ell$ are conjugate in $\G^F$, we have ${\rm Irr}_{K'}(\G^F,e_{s^\ell}^{\G^F})={\rm Irr}_{K'}(\G^F,e_s^{\G^F})$. So we have ${}^\sigma b=b$, which implies that $b\in \Z_\ell G$. Since $\Z_\ell$ $(=\Z_3)$ does not contain any primitive $4$-th root of unity, we see that the values of $\eta$ are not contained in $\Z_\ell$.

 }
\end{example}

\bigskip\noindent\textbf{Acknowledgement.} We wish to thank Professor Joseph Chuang for suggesting to consider the descent question on unipotent blocks of finite general linear groups. We also thank the referees for many very helpful comments and corrections which improved this paper a lot. 

X. Huang acknowledges support by China Postdoctoral Science Foundation (No. 2023T 160290). P. Li acknowledges support by NSFC (No. 12225108). J. Zhang acknowledges support by National Key R\&D Program of China (No. 2020YFE 0204200).



\begin{thebibliography}{1}
\expandafter\ifx\csname url\endcsname\relax
  \def\url#1{\texttt{#1}}\fi
\expandafter\ifx\csname urlprefix\endcsname\relax\def\urlprefix{URL }\fi
\expandafter\ifx\csname href\endcsname\relax
  \def\href#1#2{#2} \def\path#1{#1}\fi


\bibitem{Benson}
D. J. Benson, Representations and Cohomology I: Basic Representation theory of finite groups and associative algebras, Cambridge Studies in Advanced Mathematics, vol. {\bf30}, Cambridge University Press, 1991.

\bibitem{Boltje}
R. Boltje, The Brou\'{e} invariant of a $p$-permutation equivalence, arXiv:2007.13936v2.

\bibitem{BDR}
C. Bonnaf\'{e}, J.-F. Dat, R. Rouquier, Derived categories and Deligne--Lusztig varieties II, Ann. Math. (2) {\bf 185} (2017) 609--676.

\bibitem{Broue86}
M. Brou\'{e}, Les $\ell$-blocs des groupes ${\rm GL}(n,q)$ et ${\rm U}(n,q^2)$ et leurs structures locales, S\'{e}minaire Bourbaki Ast\'{e}risque {\bf 133--134}, 159--188, Math. France, Paris, 1986.

\bibitem{BM}
M. Brou\'{e}, J. Michel, Blocs et s\'{e}ries de Lusztig dans un groupe r\'{e}ductif fini, J. Reine Angew. Math. {\bf 395} (1989) 56--67.

\bibitem{CE}
M. Cabanes, M. Enguehard, Representation Theory of Finite Reductive Groups, New Mathematical Monogtaphs, vol. {\bf 1}., Cambridge University Press, 2004.


\bibitem{CK}
J. Chuang, R. Kessar, Symmetric groups, wreath products, Morita equivalences, and Brou\'{e}'s abelian defect group conjecture, Bull. London Math. Soc. {\bf 34} (2002) 174--184.

\bibitem{CR}
J. Chuang, R. Rouquier, Derived equivalences for symmetric groups and $\mathfrak{sl}_2$-categorification, Ann. Math. {\bf 167} (2008) 245--298.

\bibitem{Dagger}
S.W. Dagger, On the blocks of Chevalley groups, J. London Math. Soc. (2) {\bf3} (1971), 21--29.

\bibitem{DM}
F. Digne, J. Michel, Representations of Finite Groups of Lie Type, 2nd ed., London Math. Soc. Student Texts, vol. {\bf 95}, Cambridge University Press, 2020.

\bibitem{DVVcl}
O. Dudas, M. Varagnolo and E. Vasserot, Categorical actions on unipotent representations of finite classical groups, in Categorification and Higher Representation Theory, 41--104,
Contemp. Math., {\bf683}, Amer. Math. Soc., Providence, RI, 2017.

\bibitem{DVVuni}
O. Dudas, M. Varagnolo and E. Vasserot, Categorical actions on unipotent representations of finite unitary groups, Publ. Math. Inst. Hautes \'{E}tudes Sci.  {\bf129} (2019), 129--197.

\bibitem{Farrell}
N. Farrell, On the Morita Frobenius numbers of blocks offinite reductive groups, J. Algebra {\bf 471} (2017) 299--318.

\bibitem{Geck}
M. Geck, Character values, Schur indices and character sheaves, Represent. theory {\bf 7} (2003) 19--55.

\bibitem{Gr}
I. Grojnowski, Affine $\hat{\mathfrak{sl}}_p$ controls the modular representation theory of the symmetric groups and related Hecke algebras, arXiv:math/9907129.


\bibitem{H24}
X. Huang, Virtual Morita equivalences and Brauer character bijections, to appear in Arch. Math. (2024) https://doi.org/10.1007/s00013-024-02010-z..

\bibitem{HLZ}
X. Huang, P. Li, J. Zhang, The strengthened Brou\'{e} abelian defect group conjecture for ${\rm SL}(2,p^n)$ and ${\rm GL}(2,p^n)$, J. Algebra {\bf 633} (2023) 114--137.

\bibitem{IN}
I. M. Isaacs, G. Navarro, New refinements of the McKay conjecture for arbitrary finite groups, Ann. Math. (2) {\bf 156} (2002) 333--344.

\bibitem{Kessar_Linckelmann}
R. Kessar, M. Linckelmann, Descent of equivalences and character bijections, Geometric and topological aspects of the representation theory of finite groups, Springer Proc. Math. Stat., {\bf 242} (2018) 181-212.


\bibitem{Linckelmann}
M. Linckelmann, The Block Theory of Finite Group Algebras I, London Math. Soc. Student Texts, vol. {\bf 91}, Cambridge University Press, 2018.



\bibitem{Marcusonequivalences}
A. Marcus, On equivalences between blocks of group algebras: reduction to the simple components, J. Algebra {\bf184} (1996) 372--396.

\bibitem{Mi01}
H. Miyachi, Unipotent Blocks of Finite General Linear Groups in Non-defining Characteristic, PhD
thesis, Chiba university, 2001.





\bibitem{Navarro}
G. Navarro, The McKay conjecture and Galois automorphisms, Ann. Math. {\bf160} (2004) 1129--1140.

\bibitem{Rickard94}
J. Rickard, Finite group actions and \'{e}tale cohomology, Inst. Hautes \'{E}tudes Sci. Publ. Math. {\bf 80} (1994), 81--94.

\bibitem{Rickardsplendid}
J. Rickard, Splendid equivalences: derived categories and permutation modules, Proc. Lond. Math. Soc. {\bf 72} (1996) 331--358.

\bibitem{Rou}
R. Rouquier, The derived category of blocks with cyclic defect groups, in: Derived Equivalences
for Group Rings (S. K\"{o}nig, A. Zimmermann), Lecture Notes in Math. {\bf1685},
Springer Verlag, Berlin-Heidelberg, 1998, 199-220.

\bibitem{Rou2002etalas}
R. Rouquier, Complexes de cha\^{\i}nes \'{e}tales et courbes de Deligne--Lusztig, J. Algebra {\bf257} (2002) 482--508.



\bibitem{Ruh}
L. Ruhstorfer, Jordan decomposition for the Alperin--McKay conjecture, PhD thesis, Bergische Universit\"{a}t Wuppertal, 2020.

\bibitem{Ruh21}
L. Ruhstorfer, The Navarro refinement of the McKay conjecture for finite groups of Lie type in defining characteristic, J. Algebra {\bf 582} (2021) 177--205.

\bibitem{Serre}
J.P. Serre, Local fields, Graduate Texts in Mathematics {\bf 67}, 1979.

\bibitem{Turner}
W. Turner, Equivalent blocks of finite general linear groups in non-describing characteristic, J. Algebra {\bf 247} (2002) 244--267.



\bibitem{T13}
A. Turull, Strengthened Alperin--McKay conjecture for $p$-solvable groups, J. Algebra {\bf 394} (2013) 79--91.

\bibitem{WZZ}
C. Wu, K. Zhang, Y. Zhou, Blocks with defect group $\Z_{2^n}\times \Z_{2^n}\times \Z_{2^m}$, J. Algebra {\bf 510} (2018) 469--498.

\end{thebibliography}
\end{document}